\documentclass[11pt]{amsart}
\usepackage[utf8]{inputenc}
\usepackage{mathptmx} % rm & math
\usepackage{caption}
\usepackage[scaled=0.90]{helvet} % ss
\usepackage{courier} % tt

\usepackage[T1]{fontenc}
\usepackage{float}
\usepackage[%pagebackref,
hypertexnames=false, colorlinks, 
linkcolor=black,citecolor=black,urlcolor=black, linktocpage=true
]%
{hyperref} %,hypertexnames=false,colorlinks,[pagebackref]
\usepackage{amsmath}
\usepackage{wrapfig}
\usepackage{amsthm}
\usepackage{amsfonts}
\usepackage{amssymb}
\usepackage{graphicx}

\usepackage{mathrsfs}

\usepackage[margin=0.96in]{geometry}
\usepackage{color}
\usepackage{caption}

\newtheorem{lemma}{Lemma}

\newtheorem*{question}{Question}
\newcommand{\sgn}{\operatorname{sgn}}

\begin{document}

\title{A bundling problem revisited}

\author{Paata Ivanisvili}
%\thanks{PI is partially supported by the Hausdorff Institute for Mathematics, Bonn, Germany}
\address{Department of Mathematics,  Kent State University,
Kent, OH 44240, USA}
\email{ivanishvili.paata@gmail.com}

%\author{Alexander Volberg}
%\thanks{AV is partially supported by the NSF grant DMS-1265549 and by the Hausdorff Institute for Mathematics, Bonn, Germany}
%\address{Department of Mathematics,  Michigan State University, East
%Lansing, MI 48824, USA}
%\email{volberg@math.msu.edu}

\makeatletter
\@namedef{subjclassname@2010}{
  \textup{2010} Mathematics Subject Classification}
\makeatother


% 42B	Harmonic analysis in several variables
% 42B20	Singular and oscillatory integrals (Calder?on-Zygmund, etc.)
% 42B35	Function spaces arising in harmonic analysis

% 47A	General theory of linear operators
% 47A30	Norms (inequalities, more than one norm, etc.)

%{30E20, 47B37, 47B40, 30D55.}
%
% 30D55	$H^p$-classes (1980-2009)
% 30E20	Integration, integrals of Cauchy type, integral representations of analytic functions
%
% 47B   	Special classes of linear operators
% 47B37	Operators on special spaces (weighted shifts, operators on sequence spaces, etc.)
% 47B40	Spectral operators, decomposable operators, well-bounded operators, etc.


\begin{abstract}
It was conjectured by M. Glasser and S. Davison and later proved by A. Eremenko  that the certain animals should  gather close to each other in order to decrease the total heat loss. In this paper we show that it is not always true for the individual heat loss. This gives a negative answer to a question posed by A. Eremenko. 
\end{abstract}

\date{}
\maketitle

\section{Motivation, model and results}

In \cite{GD}, \cite{GD0} Glasser and Davison consider the following problem: let $B_{1}$ and $B_{2}$ are two disjoint balls in $\mathbb{R}^{3}$ with equal radii.  Let $d\geq 0$ be a distance between the  balls. Consider a harmonic function in the complement of the balls, i.e., $\mathbb{R}^{3}\setminus B_{1} \cup B_{2}$, such that $u|_{\partial B_{j}}=1$ and $\lim_{|x| \to \infty} u(x)=0$. 
Let  
\begin{align*}
Q(d) = \int_{\partial B_{1}} \frac{\partial u}{\partial n} d\sigma +\int_{\partial B_{2}} \frac{\partial u}{\partial n} d\sigma
\end{align*}
be the heat flux where $n$ is the outward unit normal vector to a sphere. Is it true that the quantity $Q(d)$ is increasing as the function of distance between the balls? 

The problem arises from the question why certain warm blooded animals like armadillos
can keep each other warm by huddling together. In this simple  model the balls $B_{1}$ and $B_{2}$ represent uniform spherical animals  in $\mathbb{R}^{3}$ with body temperature 1 and medium temperature 0.   
The harmonic function $u(x)$ represents time independent  temperature in $\mathbb{R}^{3} \setminus B_{1}\cup B_{2}$, and the quantity $Q(d)$ represents the total heat loss of both animals: the  total amount of heat given off by the animals as a function of distance between the balls. Presumably  moving  animals closer together decreases  the heat loss $Q(d)$, and it was confirmed numerically in  \cite{GD}  that the quantity $Q(d)$ is increasing. However no mathematical proof was given until A.~Eremenko \cite{Er} gave a rigorous proof in more general setting. In Eremenko's argument it was noticed that 
\begin{align}\label{cap}
\frac{4\pi}{Q}=\inf\left\{ I(\mu) \; : \mathrm{supp} \, \mu \subset \partial B_{1} \cup \partial B_{2}, \; \mu(\partial B_{1} \cup \partial B_{2})=1, \; \mu \geq 0\right\}
\end{align}
where 
\begin{align*}
I(\mu) = \int \int \frac{d\mu(x) d\mu(y)}{|x-y|}.
\end{align*}
The monotonicity of $Q$ follows from the fact that if $\varphi : \mathbb{R}^{3} \to \mathbb{R}^{3}$ is continuous,  one-to-one, and $|\varphi(x) - \varphi(y)| \geq  |x-y|$ then $I(\varphi_{*}\mu) \leq I(\mu)$ where $\varphi_{*}\mu(A) = \mu(\varphi^{-1}(A))$ for any Borel measurable $A \subset \mathbb{R}^{3}$. 

 Notice that each individual animal $B_{j}$ feels only his own heat loss $Q_{j}(d)$
\begin{align*}
Q_{j}(d) = \int_{\partial B_{j}} \frac{\partial u}{\partial n} d\sigma,
\end{align*}
but not the total $Q=Q_{1}+Q_{2}$. Therefore the behavior of the animals we have discussed could be driven by individual feelings but not the abstract ``common goal''. In case of equal balls  the individual heat loss $Q_{j}$ is monotonically increasing because $Q_{1}=Q_{2} = \frac{Q}{2}$.  It is natural to think that the individual heat loss $Q_{j}$ is monotonically increasing  for the balls of different radii. In \cite{Er, Er2} the following question was asked:
\begin{question}
Are the quantities $Q_{j}(d),$ $j=1,2$ monotonically increasing if the balls have different radii? 
\end{question}

In Section~\ref{computation} we will show that if $r_{j}$ denotes the radius of the ball $B_{j}$ for $j=1,2$, then the of heat loss $Q_{1}$ is not monotonically increasing provided that $\frac{r_{1}}{r_{2}}>\ell $ where $\ell\approx 1.95$ is the positive  solution of the equation
\begin{align}\label{shedegi}
2(1+x^{3})\left(\gamma+\psi\left(\frac{1}{1+x}\right)\right)+x^{2}+x+2(x^{2}-x)\psi'\left(\frac{1}{1+x}\right)=0.
\end{align}
In the above notation  $\gamma$ is the Euler's constant, and $\psi$ is the digamma function. 

In other words, if armadillo A is at least twice as big as armadillo B, then A should keep some nonzero distance from B in order to minimize the heat loss while B should try to be as close as possible to A. 

In Section~\ref{computation} we obtain the following asymptotic expression for $Q_{1}(d)$ which implies the conclusion (\ref{shedegi}). Let $r_{1}>r_{2}>0$, $d\geq 0$,  and let the temperature of the balls be constant and equal to $T_{0}>0$. Then 
\begin{align*}
 &\frac{Q_{1}(d)}{4 \pi T_{0} r_{1}} = - \frac{r_{2}\left(\gamma+ \psi \left( \frac{r_{2}}{r_{1}+r_{2}}\right)\right)}{r_{1}+r_{2}}-\\
  &\frac{d}{6r_{1}(r_{1}+r_{2})^{3}}\left[2(r_{1}^{3}+r_{2}^{3})\left(\gamma+\psi\left(\frac{r_{2}}{r_{1}+r_{2}}\right)\right) +r_{1}^{2}r_{2}+r_{2}^{2}r_{1}+2(r_{1}^{2}r_{2}-r_{2}^{2}r_{1})\psi'\left(\frac{r_{2}}{r_{1}+r_{2}} \right) \right]+o(d).
\end{align*}
as $d \to 0$. 

\section{Two balls of unequal radii}\label{computation}
We consider a bispherical coordinate system
\begin{align*}
(x,y,z)=\left( \frac{a\sin \eta \cos \phi}{\cosh \mu - \cos \eta}, \frac{a \sin \eta \sin \phi}{\cosh \mu \cos \eta}, \frac{a \sin \mu}{\cosh \mu - \cos \eta}\right),
\end{align*}
where  $-\infty\leq \mu \leq \infty$, $0 \leq \eta \leq \pi$, $0 \leq \varphi \leq 2\pi$ and $a>0$ so that the foci $F_{1}$ and $F_{2}$ coincide with the centers of  $B_{1}$ and $B_{2}$. A general solution of the Laplace equation in  bispherical coordinate system  (under the assumptions that the solution does not depend on $\varphi$ which is true in our case) is given by the expression (see page 1298, \cite{MF})
\begin{align*}
T(\mu, \eta)=\sqrt{\cosh \mu - \cos \eta} \sum_{n=0}^{\infty}\left[A_{n} e^{(n+\frac{1}{2})\mu}+B_{n} e^{-(n+\frac{1}{2})\mu} \right]P_{n}(\cos \eta), 
\end{align*}
where $P_{n}$ are Legandre polynomials. Let $r_{1}>0$ and $r_{2}>0$ be the radii of $B_{1}$ and $B_{2}$ correspondingly. The corresponding values of the $\mu$ coordinate  $\mu=\mu_{1}>0$ on $\partial B_{1}$ and  $\mu=\mu_{2}<0$ on $\partial B_{2}$ will be determined by $ r_{j} = \frac{a}{|\sinh \mu_{j}|}$ for $j=1,2$. Notice also that the distance $d$ between the balls can be obtained as follows  
\begin{align*}
d+r_{1}+r_{2}=a(\coth \mu_{1}-\coth \mu_{2}) = r_{1} \cosh \mu_{1}+r_{2}\cosh \mu_{2}. 
\end{align*}

Assume that the temperature on $\partial B_{j}$ is a constant $T_{j}$, $j=1,2$. Then by using the generating function for the Legandre polynomials 
\begin{align}\label{gen}
\frac{1}{\sqrt{\cosh \mu - \cos \eta}} = \sqrt{2} \sum_{n=0}^{\infty} e^{-(n+\frac{1}{2})|\mu|}P_{n}(\cos \eta),
\end{align}
the boundary condition $T(\mu_{j}, \eta)=T_{j}$ implies 
\begin{align}\label{boun}
A_{n} e^{(n+\frac{1}{2})\mu_{j}}+B_{n} e^{-(n+\frac{1}{2})\mu_{j}}=T_{j} \sqrt{2} e^{-(n+\frac{1}{2})|\mu_{j}|}, \quad j=1,2 \quad n \geq 0. 
\end{align}

Let us compute the surface element $d\sigma$ in the bispherical coordinates:
\begin{align*}
| (x,y,z)_{\eta}\times (x,y,z)_{\phi} | = \frac{a^{2} \sin \eta}{(\cosh \mu - \cos \eta)^{2}}.
\end{align*}
The heat loss $Q_{1}$ can be computed as follows  
\begin{align*}
Q_{1} = \iint_{\partial B_{1}} \left.\frac{\partial T}{\partial n}\right|_{\mu=\mu_{1}} \frac{a^{2}\sin \eta d\eta d\phi}{(\cosh \mu_{1}-\cos \eta)^{2}}.
\end{align*}
Notice that on $\partial B_{1}$ we have
\begin{align}
&\frac{\partial T}{\partial n} = \left(\frac{\cosh \mu_{1} - \cos \eta}{a}\right) \frac{\partial T}{\partial \mu} =\nonumber \\
&\frac{\sinh \mu_{1}}{2a}\sqrt{\cosh \mu_{1}-\cos \eta}\sum_{n=0}^{\infty}\left[A_{n} e^{(n+\frac{1}{2})\mu_{1}}+B_{n} e^{-(n+\frac{1}{2})\mu_{1}} \right]P_{n}(\cos \eta)+\nonumber \\
&\frac{(\cosh \mu_{1} -\cos \eta)^{3/2}}{a}\sum_{n=0}^{\infty}\left[(n+\frac{1}{2})A_{n}e^{(n+\frac{1}{2})\mu_{1}}-(n+\frac{1}{2})B_{n}e^{-(n+\frac{1}{2})\mu_{1}} \right]P_{n}(\cos \eta)=\nonumber \\
&\frac{\sinh \mu_{1}}{2a}T(\mu_{1}, \eta)+ \frac{(\cosh \mu_{1} - \cos \eta)^{3/2}}{a}\sum_{n=0}^{\infty}(n+\frac{1}{2})\left[ A_{n}e^{(n+\frac{1}{2})\mu_{1}}-B_{n}e^{-(n+\frac{1}{2})\mu_{1}}\right]P_{n}(\cos \eta).\label{last}
\end{align}
The boundary conditions (\ref{boun}) imply that 
\begin{align*}
&A_{n}e^{(n+\frac{1}{2})\mu_{1}}-B_{n} e^{-(n+\frac{1}{2})\mu_{1}} = \sqrt{2} \frac{e^{-(n+\frac{1}{2})\mu_{1}}T_{1}+e^{-(n+\frac{1}{2})(3\mu_{1}-2\mu_{2})}T_{1}-2T_{2}e^{-(n+\frac{1}{2})(\mu_{1}-2\mu_{2})}}{1-e^{-(n+\frac{1}{2})(2\mu_{1}-2\mu_{2})}} =\sum_{k=0}^{\infty}\\
&\left(T_{1}\sqrt{2}e^{-(n+\frac{1}{2})(\mu_{1}+k(2\mu_{1}-2\mu_{2}))}+T_{1}\sqrt{2}e^{-(n+\frac{1}{2})(3\mu_{1}-2\mu_{2}+k(2\mu_{1}-2\mu_{2}))}-T_{2}2\sqrt{2}e^{-(n+\frac{1}{2})(\mu_{1}-2\mu_{2}+k(2\mu_{1}-2\mu_{2}))}\right)=\\
&T_{1}\sqrt{2}e^{-(n+\frac{1}{2})\mu_{1}}+2\sqrt{2}\sum_{k=0}^{\infty}\left(T_{1}e^{-(n+\frac{1}{2})(3\mu_{1}-2\mu_{2}+k(2\mu_{1}-2\mu_{2}))}-T_{2}e^{-(n+\frac{1}{2})(\mu_{1}-2\mu_{2}+k(2\mu_{1}-2\mu_{2}))} \right).
\end{align*}

Thus (\ref{last}) takes the following form 
\begin{align*}
&\frac{\partial T}{\partial n} = \frac{\sinh \mu_{1}}{2a}T(\mu_{1}, \eta)+ \frac{(\cosh \mu_{1} - \cos \eta)^{3/2}}{a} T_{1} \sqrt{2}\sum_{n=0}^{\infty} (n+\frac{1}{2}) e^{-\left(n+\frac{1}{2}\right)\mu_{1}}P_{n}(\cos \eta)+\\
&\frac{(\cosh \mu_{1} - \cos \eta)^{3/2}}{a}\cdot  2\sqrt{2} T_{1} \sum_{n=0}^{\infty}\sum_{k=0}^{\infty} (n+\frac{1}{2}) e^{-(n+\frac{1}{2})(3\mu_{1}-2\mu_{2}+k(2\mu_{1}-2\mu_{2}))} P_{n}(\cos \eta) \\
&-\frac{(\cosh \mu_{1} - \cos \eta)^{3/2}}{a}\cdot  2\sqrt{2} T_{2} \sum_{n=0}^{\infty}\sum_{k=0}^{\infty} (n+\frac{1}{2}) e^{-(n+\frac{1}{2})(\mu_{1}-2\mu_{2}+k(2\mu_{1}-2\mu_{2}))}  P_{n}(\cos \eta).
\end{align*}

If we differentiate (\ref{gen}) with respect to  $\mu$ we obtain   the following identity
\begin{align}\label{tozhd}
\frac{1}{2}\sinh \mu (\cosh \mu - \cos \eta)^{-3/2}=\sqrt{2} \sum_{n=0}^{\infty}(n+\frac{1}{2})e^{-(n+\frac{1}{2})\mu}P_{n}(\cos \eta) \quad \text{for} \quad \mu>0.
\end{align}
Using (\ref{tozhd}) we  further simplify the expression for $\frac{\partial T}{\partial n}$:  
\begin{align*}
&\frac{\partial T}{\partial n}=\frac{\sinh \mu_{1}}{2a}T(\mu_{1}, \eta) +  \frac{\sinh \mu_{1}}{2a}T_{1}+\\
& \sum_{k=0}^{\infty}\frac{T_{1}}{a}\sinh (3\mu_{1}-2\mu_{2}+k(2\mu_{1}-2\mu_{2}))\cdot \left( \frac{\cosh \mu_{1} - \cos \eta}{\cosh(3\mu_{1} - 2\mu_{2}+k(2\mu_{1}-2\mu_{2}))-\cos \eta} \right)^{3/2}\\
&- \sum_{k=0}^{\infty}\frac{T_{2}}{a}\sinh(\mu_{1}-2\mu_{2}+k(2\mu_{1}-2\mu_{2})) \left(\frac{\cosh \mu_{1} - \cos \eta}{\cosh(\mu_{1}-2\mu_{2}+k(2\mu_{1}-2\mu_{2}))-\cos \eta}\right)^{3/2}
\end{align*}
We notice that $T(\mu_{1},\eta) = T_{1}$. Then $\frac{\sinh \mu_{1}}{2a}T(\mu_{1}, \eta) +  \frac{\sinh \mu_{1}}{2a}T_{1} = \frac{\sinh \mu_{1}}{a} T_{1}$ and  finally we have  
\begin{align*}
&Q_{1}= \int_{0}^{\pi} \left[ 2 \pi a T_{1} \sinh \mu_{1} \frac{\sin \eta }{(\cosh \mu_{1} - \cos \eta)^{2}} + \sum_{k=0}^{\infty}\right.\\
&\left( 2 \pi a T_{1} \sinh (3 \mu_{1} - 2\mu_{2} + k(2\mu_{1} -2\mu_{2})) \cdot \frac{(\cosh \mu_{1} - \cos \eta)^{-1/2}\sin \eta }{(\cosh(3\mu_{1}-2\mu_{2}+k(2\mu_{1}-2\mu_{2}))-\cos \eta)^{3/2}}-\right.\\
&\left. \left. 2 \pi a T_{2} \sinh (\mu_{1}-2\mu_{2}+k(2\mu_{1}-2\mu_{2}))\cdot \frac{(\cosh \mu_{1} -\cos \eta)^{-1/2}\sin \eta}{(\cosh(\mu_{1} -2\mu_{2}+k(2\mu_{1}-2\mu_{2}))-\cos \eta)^{3/2}} \;  \right)\right]  d\eta .
\end{align*}
By substituting  $\cos \eta = x$ we obtain 
\begin{align*}
&Q_{1} = \int_{-1}^{1}  \left[ 2 \pi a T_{1} \sinh \mu_{1} \frac{1}{(\cosh \mu_{1}-x)^{2}}+\sum_{k=0}^{\infty}\right. \\
&\left(2 \pi a T_{1} \sinh (3\mu_{1} - 2\mu_{2} + k(2\mu_{1} -2\mu_{2})) \cdot \frac{(\cosh \mu_{1} - x)^{-1/2}}{(\cosh (3\mu_{1} - 2\mu_{2} +k(2\mu_{1} -2\mu_{2})) -x)^{3/2}}-\right.\\
&\left. \left.  2 \pi a T_{2} \sinh (\mu_{1} - 2\mu_{2} + k(2\mu_{1} -2\mu_{2}))\cdot \frac{(\cosh \mu_{1} -x)^{-1/2}}{(\cosh (\mu_{1} - 2\mu_{2} +k(2\mu_{1} - 2\mu_{2}))-x)^{3/2}}\; \right) \right]dx =\\
&I_{1} + \sum_{k=0}^{\infty}(I_{2,k}+I_{3,k}). 
\end{align*}
Let us calculate each term separately. We remind that  $\frac{\sinh \mu_{1}}{2a} T_{1}=\frac{T_{1}}{2r_{1}}$. Therefore we have 
\begin{align*}
I_{1}=\frac{4 \pi a T_{1} \sinh \mu_{1}}{\sinh^{2} \mu_{1}}=4 \pi T_{1} r_{1}.
\end{align*}
We notice the following subtle identity 
\begin{align}\label{ident}
\sinh(B) \int_{-1}^{1} \frac{(\cosh(A) -x)^{-1/2} dx}{(\cosh (B) -x)^{3/2}}=\frac{2}{\sinh\left( \frac{A+B}{2}\right)}
\end{align}
for all real numbers  $A$ and $B$ whenever the both sides of (\ref{ident}) make sense. 
Then 
\begin{align*}
&I_{2,k}= \frac{4 \pi a T_{1}}{\sinh(2\mu_{1}-\mu_{2}+k(\mu_{1}-\mu_{2}))},\\
&I_{3,k}=-\frac{4 \pi a T_{2}}{\sinh(\mu_{1}-\mu_{2}+k(\mu_{1}-\mu_{2}))}.
\end{align*}
Taking into account that $r_{1} = \frac{a}{\sin \mu_{1}}$ we obtain 
\begin{align*}
&Q_{1} = 4 \pi T_{1} r_{1} + 4 \pi a \sum_{k=0}^{\infty} \frac{T_{1}}{\sinh(2\mu_{1}-\mu_{2}+k(\mu_{1}-\mu_{2}))} - \frac{T_{2}}{\sinh(\mu_{1}-\mu_{2}+k(\mu_{1}-\mu_{2}))}=\\
&4 \pi T_{1} r_{1} + 4\pi r_{1} \sinh(\mu_{1}) \sum_{k=1}^{\infty} \frac{T_{1}}{\sinh(\mu_{1}+k(\mu_{1}-\mu_{2}))}-\frac{T_{2}}{\sinh(k(\mu_{1}-\mu_{2}))}.
\end{align*}
Further we consider the case when $T_{1}=T_{2}=T_{0}>0$. In order to investigate the monotonicity of $Q_{1}$ with respect to $d$, it is enough to investigate the monotonicity of the following function 
\begin{align*}
f(d) := \sinh(\mu_{1}) \sum_{k=1}^{\infty} \frac{1}{\sinh(\mu_{1}+k(\mu_{1}-\mu_{2}))}-\frac{1}{\sinh(k(\mu_{1}-\mu_{2}))}.
\end{align*} 
We notice that 
\begin{align*}
\cosh \mu_{1} = \frac{(d+r_{1}+r_{2})^{2}+r_{1}^{2}-r_{2}^{2}}{2(d+r_{1}+r_{2})r_{1}} \quad \text{and} \quad \cosh \mu_{2} = \frac{(d+r_{1}+r_{2})^{2}+r_{2}^{2}-r_{1}^{2}}{2(d+r_{1}+r_{2})r_{2}}.
\end{align*}
Let $x=e^{\mu_{1}}>1$ and $y = e^{-\mu_{2}} >1$. Then $f(d)$ takes the form 
\begin{align}\label{bolo}
f(d) = \left(x-\frac{1}{x} \right) \sum_{k=1}^{\infty} \left( \frac{x(xy)^{k}}{x^{2}(xy)^{2k}-1} - \frac{(xy)^{k}}{(xy)^{2k}-1}\right)
\end{align}
where 
\begin{align}
&x =1+ \frac{d(d+2r_{2})+\sqrt{(d+2r_{1}+2r_{2})(d+2r_{2})(d+2r_{1})d}}{2(d+r_{1}+r_{2})r_{1}}; \label{ixi}\\
&y=1+\frac{d(d+2r_{1})+\sqrt{(d+2r_{1}+2r_{2})(d+2r_{2})(d+2r_{1})d}}{2(d+r_{1}+r_{2})r_{2}}\label{yi}. 
\end{align}
By using the identity $\frac{1}{t-1} = \frac{1}{t}\sum_{k=0}^{\infty}t^{-k}$ two times for the terms inside the summation (\ref{bolo}), and  by Fubini's theorem  the expression for $f(d)$ can be simplified as follows 
\begin{align*}
f(d) = \left( x-\frac{1}{x} \right) \sum_{j=0}^{\infty}  \frac{x^{-2j-1}-1}{(xy)^{2j+1}-1}.
\end{align*}
If $r_{1}=r_{2}$ then, as we already mentioned in (\ref{cap}), it is known that $f$ is monotonically increasing (see also a proof in \cite{fp} without resorting to (\ref{cap})). Therefore it is enough to study the sign of $\lim_{d \to 0+}f'(d)$ for different radii $r_{1}$ and $r_{2}$. Further we assume that $r_{1}>r_{2}$. 
  Let  $z=xy$. 
 Set  
\begin{align*}
f(d) = \left( x -\frac{1}{x}\right) \sum_{j=0}^{\infty}\frac{x^{-2j-1}-1}{z^{2j+1}-1} = \left( x -\frac{1}{x}\right) \sum_{k=0}^{\infty} g(k),
\end{align*}
where $g(s) = \frac{x^{-2s-1}-1}{z^{2s+1}-1}$. Since $x, z >1$ it is easy to see that  $g(s) \in C^{\infty}([0,\infty))$ and all its derivatives tend to zero as $k \to \infty$. Therefore By Euler--Maclaurin formula we have 
\begin{align}\label{euler}
\left(x-\frac{1}{x} \right)\sum_{k=0}^{\infty}g(k)=\left(x-\frac{1}{x} \right)\left( \int_{0}^{\infty} g(s) ds +\frac{g(0)}{2}-\frac{g'(0)}{12}\right) -\left(x-\frac{1}{x} \right) \int_{0}^{\infty} \frac{B_{2}(\{1-s\})}{2} g^{(2)}(s)ds 
\end{align} 
where $B_{2}(x)=x^{2}-x+\frac{1}{2}$ is the  Bernoulli polynomial, and $\{x\}$ represents the fractional part of $x$. 

We will compute the asymptotic behavior  of each term in (\ref{euler}) separately as $t \to 0$. First notice that (\ref{ixi}) implies
\begin{align}\label{pirveli}
x-\frac{1}{x}=2\sqrt{\frac{2r_{2}}{(r_{1}+r_{2})r_{1}}} \, d^{1/2}+\frac{1}{\sqrt{2}} \frac{r_{1}^{2}-r_{1}r_{2}+r_{2}^{2}}{(r_{1}+r_{2})^{3/2}\sqrt{r_{2}} r_{1}^{3/2}} d^{3/2}+O(d^{2}).
\end{align} 
We have 
\begin{align*}
\int_{0}^{\infty}g(t)dt = \int_{0}^{\infty} \frac{x^{-1} e^{-2t \ln x} -1}{ze^{2 t \ln z}-1}dt = \frac{1}{ \ln z} \int_{0}^{\infty} \frac{x^{-1} e^{-s\frac{\ln x}{\ln z}}-1}{ze^{s}-1}ds. 
\end{align*}
We notice that 
\begin{align*}
\int_{0}^{\infty} \frac{1}{ze^{s}-1}ds = \ln \left(\frac{z}{z-1}\right).
\end{align*}
Therefore 
\begin{align*}
\frac{1}{ \ln z} \int_{0}^{\infty} \frac{x^{-1} e^{-s\frac{\ln x}{\ln z}}-1}{ze^{s}-1}ds = -\frac{\ln \left( \frac{z}{z-1}\right)}{2\ln z} + \frac{1}{2x \ln z} \int_{0}^{\infty} \frac{e^{-s\frac{\ln x}{\ln z}}}{ze^{s}-1}ds.
\end{align*}
Set $\alpha = \frac{\ln x}{\ln z}$.  Then 
\begin{align*}
&\int_{0}^{\infty} \frac{e^{-s\alpha }}{ze^{s}-1}ds = \int_{1}^{\infty}y^{-\alpha}\left[-\frac{1}{y} +\frac{z}{zy-1}\right]dy = -\frac{1}{\alpha} + z \int_{1}^{\infty} \frac{y^{-\alpha}}{ zy-1}dy=\\
&-\frac{1}{\alpha}+ x\int_{0}^{1/z} \frac{s^{\alpha-1}}{1-s} ds = \frac{1}{z} \sum_{k=0}^{\infty} \frac{(z^{-1})^{k}}{\alpha+1+k}=\frac{z^{-1}}{1+\alpha}\;  {}_{2}F_{1}(1,1+\alpha;2+\alpha;z^{-1}),
\end{align*}
where 
\begin{align*}
{}_{2}F_{1}(a,b; c; z) = \sum_{n=0}^{\infty} \frac{(a)_{n} (b)_{n}}{(c)_{n}} \frac{z^{n}}{n!}
\end{align*}
 denotes the hypergeometric function where $(a)_{n} = a(a+1)\cdots (a+n-1)$ if $n \geq 1$ and $(a)_{0}=1$. On the other hand it is known that (see~\cite{nume})
 \begin{align*}
 {}_{2}F_{1}(a,b; a+b; v)= \frac{\Gamma(a+b)}{\Gamma(a) \Gamma(b)} \left(  \sum_{k=0}^{\infty} \frac{(a)_{k} (b)_{k}}{k!^{2}} (-\ln (1-v)+2\psi(k+1)-\psi(a+k)-\psi(b+k))(1-v)^{k}\right)  
 \end{align*}
 for all $0<1-v<1$, where $\psi$ is the digamma function.  Therefore we obtain 
 \begin{align*}
& \frac{1}{2x \ln z} \int_{0}^{\infty} \frac{e^{-s \frac{\ln x}{\ln z}}}{ze^{s}-1}ds = \frac{1}{2xz \ln z} \frac{{}_{2}F_{1}(1,1+\alpha;2+\alpha;z^{-1})}{1+\alpha} =\\
 &\frac{1}{2x z \ln z} \sum_{k=0}^{\infty}\left( \frac{(1+\alpha)_{k}}{k!}(-\ln(1-z^{-1})+\psi(k+1)-\psi(1+k+\alpha))(1-z^{-1})^{k}\right).
 \end{align*}
 Note that  when $d \to 0+$ we have   $\alpha  = \frac{r_{2}}{r_{1}+r_{2}}+O(\sqrt{d})$, $1-z^{-1}=\sqrt{\frac{2(r_{1}+r_{2})}{r_{1} r_{2}}} \; d^{1/2}+O(d)$ and  $\ln (1-z^{-1})  (1-z^{-1})^{3} = O(d^{3/2}\ln d)$. It is known that $\psi(t) = \ln(t) + O(1/t)$ for $t \to \infty$. Therefore if $d$  is sufficiently small we have $\alpha <1$ and thus for $k\geq 3$ we obtain
 \begin{align*}
\left| \left( \frac{(1+\alpha)_{k}}{k!}(-\ln(1-z^{-1})+\psi(k+1)-\psi(1+k+\alpha))(1-z^{-1})^{k}\right) \right| \leq 10 k b(r_{1},r_{2})| \ln d | (c(r_{1},r_{2}) d^{1/2})^{k},
 \end{align*}
 where $b(r_{1},r_{2})>0$ and $c(r_{1},r_{2})>0$ are some finite numbers depending on $r_{1}$ and $r_{2}$.
 Therefore for sufficiently small $d$ we have 
 \begin{align*}
&\left|  \sum_{k=3}^{\infty}\left( \frac{(1+\alpha)_{k}}{k!}(-\ln(1-z^{-1})+\psi(k+1)-\psi(1+k+\alpha))(1-z^{-1})^{k}\right) \right| \leq \\
&10 b(r_{1},r_{2}) | \ln d | \sum_{k=3}^{\infty} k  (c(r_{1},r_{2}) d^{1/2})^{k} \leq 100 b(r_{1},r_{2}) | \ln d | c(r_{1},r_{2})^{3} d^{3/2}= O(d^{3/2} \ln d).
 \end{align*}
 
 We obtain 
 \begin{align*}
 & \sum_{k=0}^{\infty}\left( \frac{(1+\alpha)_{k}}{k!}(-\ln(1-z^{-1})+\psi(k+1)-\psi(1+k+\alpha))(1-z^{-1})^{k}\right) = \\
 & \sum_{k=0}^{2}\left( \frac{(1+\alpha)_{k}}{k!}(-\ln(1-z^{-1})+\psi(k+1)-\psi(1+k+\alpha))(1-z^{-1})^{k}\right) + O(d^{3/2}\ln d). 
 \end{align*}
 Notice that $ (x-\frac{1}{x}) \frac{1}{2 xz \ln z} = \frac{r_{2}}{r_{1}+r_{2}}+O(\sqrt{d})$. Thus we obtain 
 \begin{align*}
 &\left(x- \frac{1}{x}\right) \left[ \int_{0}^{\infty} g(s)ds + \frac{g(0)}{2}-\frac{g'(0)}{12}\right] =\left(x-\frac{1}{x}\right)\left[\frac{g(0)}{2}-\frac{g'(0)}{12} - \frac{\ln\left(\frac{z}{z-1}\right)}{2\ln z} + \right.\\
 &\left. \frac{1}{2x z \ln z} \left( O(d^{3/2} \ln d)+\sum_{k=0}^{2}\left( \frac{(1+\alpha)_{k}}{k!}(-\ln(1-z^{-1})+\psi(k+1)-\psi(1+k+\alpha))(1-z^{-1})^{k}\right)   \right)\right]=\\
 &\left(x-\frac{1}{x} \right) \left[ \frac{x^{-1}-1}{2(z-1)} +\frac{z\ln x - \ln x + z \ln z - z x\ln z}{6x(z-1)^{2}}- \frac{\ln\left(\frac{z}{z-1}\right)}{2\ln z} \right.+\\
 &\left. \frac{1}{2x z \ln z} \sum_{k=0}^{2} \frac{(1+\alpha)_{k}}{k!}(-\ln(1-z^{-1})+\psi(k+1)-\psi(1+k+\alpha))(1-z^{-1})^{k}\right]   +O(d^{3/2} \ln d)
 %\\
 %&\frac{1}{2xz \ln z} \left(x-\frac{1}{x} \right) \left[ \frac{(1-x)z\ln z}{z-1}+ \frac{(z\ln x - \ln x + z \ln z - z x\ln z)z\ln z}{3(z-1)^{2}} - xz\ln\left(\frac{z}{z-1}\right)\right. +\\
 %&\left(1+(1-z^{-1})(1+\alpha) + \frac{(1+\alpha)(2+\alpha)(1-z^{-1})^{2}}{2}\right)\ln\left(\frac{z}{z-1} \right)+ \psi(1)-\psi(1+\alpha)+\\
 %&\left. (1+\alpha)(1-z^{-1})(\psi(2)-\psi(2+\alpha))+\frac{(1+\alpha)(2+\alpha)(1-z^{-1})^{2}}{2} (\psi(3)-\psi(3+\alpha))\right] + O(d^{3/2} \ln d)
 \end{align*}
 We note that $\psi(1)=-\gamma$, then one can check that when $d \to 0$,  after some routine  computations, using (\ref{ixi}), (\ref{yi})  and  the identity $\psi(x+1)=\psi(x)+\frac{1}{x}$ several times, the above expression takes the following form
 \begin{align*}
  &-\frac{r_{2}\left(\gamma+ \psi \left( \frac{2r_{2}+r_{1}}{r_{1}+r_{2}}\right)\right)}{r_{1}+r_{2}}-\\
  &\frac{d}{6r_{1}(r_{1}+r_{2})^{3}}\left[(r_{1}^{3}+r_{2}^{3})\left(2\gamma+2\psi\left(\frac{r_{2}}{r_{1}+r_{2}}\right)\right) +r_{1}^{2}r_{2}+r_{2}^{2}r_{1}+2(r_{1}^{2}r_{2}-r_{2}^{2}r_{1})\psi'\left(\frac{r_{2}}{r_{1}+r_{2}} \right) \right]+o(d).
 \end{align*}

 We are left with showing that the term $\left(x-\frac{1}{x} \right)\int_{0}^{\infty} \frac{B_{2}(\{ 1-s\})}{2} g^{(2)}(s) ds$ in (\ref{euler}) behaves as $o(d)$ for sufficiently small $d$. Since $x-\frac{1}{x}=O(\sqrt{d})$ it is enough to show that $\int_{0}^{\infty} \frac{B_{2}(\{ 1-s\})}{2} g^{(2)}(s) ds = o(\sqrt{d})$. We have
 \begin{align}\label{ragaca}
 \left| \int_{0}^{\infty}B_{2}( \{1-s \}) g^{(2)}(s) ds\right|  = \frac{1}{3}\left|\int_{0}^{\infty} B_{3}( \{1-s \}) g^{(3)}(s) \right| \leq \int_{0}^{\infty} |g^{(3)}| ds 
 \end{align}
 where $B_{3}(x) = x^{3} - \frac{3}{2}x^{2} + \frac{1}{2}x$ is Bernoulli polynomial.
 Consider a function $f(s) = \frac{e^{sp}-1}{e^{s}-1}$ where $p=-\frac{\ln x}{\ln xy}$. Clearly $f((2s+1)\ln xy)=g(s)$. We need the following technical lemma:
 \begin{lemma}
 Let $-1 \leq p \leq 0$. Then $f^{(3)}(s) \geq 0$ for all $s \geq 0$. 
 \end{lemma}
 
 Before we proceed to the proof of the lemma  we will show how the desired estimate follows from the lemma. 
 First notice that $p = - \frac{\ln x}{\ln xy} \in [-1,0]$ because $x,y > 1$. Therefore $\sgn( g^{(3)}) = \sgn(f^{(3)}) >0$  for all  $ s >0$,   and we have 
 \begin{align*}
 \int_{0}^{\infty} | g^{(3)}| ds = -g''(0).
 \end{align*}
 After some straightforward computations one can show that $g''(0) = O(d)$ as $d \to 0$. We will omit the details of the unnecessary computations. 
 
 Finally, we obtain that 
 \begin{align*}
 &\frac{Q_{1}(d)}{4 \pi T_{0} r_{1}} = 1- \frac{r_{2}\left(\gamma+ \psi \left( \frac{2r_{2}+r_{1}}{r_{1}+r_{2}}\right)\right)}{r_{1}+r_{2}}-\\
  &\frac{d}{6r_{1}(r_{1}+r_{2})^{3}}\left[(r_{1}^{3}+r_{2}^{3})\left(2\gamma+2\psi\left(\frac{r_{2}}{r_{1}+r_{2}}\right)\right) +r_{1}^{2}r_{2}+r_{2}^{2}r_{1}+2(r_{1}^{2}r_{2}-r_{2}^{2}r_{1})\psi'\left(\frac{r_{2}}{r_{1}+r_{2}} \right) \right]+o(d)=\\
  &-\frac{r_{2}\left(\gamma +\psi\left( \frac{r_{2}}{r_{1}+r_{2}}\right) \right)}{r_{1}+r_{2}}-\\
&\frac{d}{6r_{1}(r_{1}+r_{2})^{3}}\left[(r_{1}^{3}+r_{2}^{3})\left(2\gamma+2\psi\left(\frac{r_{2}}{r_{1}+r_{2}}\right)\right) +r_{1}^{2}r_{2}+r_{2}^{2}r_{1}+2(r_{1}^{2}r_{2}-r_{2}^{2}r_{1})\psi'\left(\frac{r_{2}}{r_{1}+r_{2}} \right) \right]+o(d)\\
 \end{align*}
 where $T_{0}=T_{1}=T_{2}$ is the temperature of the balls.

 It remains to prove the  technical lemma.
 \begin{proof}
 Notice that 
 \begin{align}
 &f^{(3)}(s) =\nonumber \\
 &\frac{e^{s(p+3)}(p-1)^{3}+e^{s(p+2)}(-3p^{3}+6p^{2}-4)+e^{s(p+1)}(3p^{3}-3p^{2}-3p-1)-p^{3}e^{sp}+e^{3s}+4e^{2s}+e^{s}}{(e^{s}-1)^{4}} \label{num}
 \end{align}
 It is enough to  show that the coefficient $a_{k}$ of $\frac{s^{k}}{k!}$ for $k \geq 0$  of  numerator in (\ref{num}) is nonnegative. Indeed 
 \begin{align*}
 &a_{k}=(p+3)^{k}(p-1)^{3}+(p+2)^{k}(-3p^{3}+6p^{2}-4)+(p+1)^{k}(3p^{3}-3p^{2}-3p-1)-p^{k+3}+3^{k}+4\cdot 2^{k}+1=\\
 &J_{1,k}(p)+J_{2,k}(p)+J_{3,k}(p),
 \end{align*}
 where 
 \begin{align*}
 &J_{1,k}=(p+1)^{k}(3p^{3}-3p^{2}-3p-1)+1;\\
 &J_{2,k}=(p+2)^{k}(-3p^{3}+6p^{2})-p^{k+3}\\
 &J_{3,k}=(p+3)^{k} (p-1)^{3}-4 (p+2)^{k}+3^{k}+4\cdot 2^{k};
 \end{align*}
 We notice that $a_{0}=a_{1}=a_{2}=a_{3}=0$, $a_{4} = 6p^{2}(p-1)^{2}\geq 0$ and  $a_{5}=4p(p-1)(6p^{2}(p+1)-14p+1)\geq 0$ because of the assumptions on $p$. It is enough to show that $J_{1,k}, J_{2,k}$ and $J_{3,k}$ are nonnegative for all $k \geq 6$. Indeed, since $p \in [-1,0]$ we have 
 \begin{align*}
 J_{1,k}(p)=(p+1)^{k}(3p^{3}-3p^{2}-3p-1)+1 \geq (p+1)(3p^{3}-3p^{2}-3p-1)+1 =p(3p^{3}-6p-4)\geq 0
 \end{align*} 
 For $J_{2,k}$ if $k$ is even then $-p^{k+3}\geq 0$ and there is nothing to prove. Assume that $k=2m+1$ where $m\geq 3$. Then 
 \begin{align*}
 J_{2,k}=(p+2)^{2m+1}(-3p^{3}+6p^{2})-p^{2m+4}\geq (p+2)(-3p^{3}+6p^{2})-p^{4} = 4p^{2}(3-p^{2}) \geq 0
 \end{align*}
It remains to show that  $J_{3,k}(p) \geq 0$. We have $J_{3,k}(0)=0$, and 
 \begin{align*}
 & J'_{3,k}(p) = k(p+3)^{k-1}(p-1)^{3}+3(p+3)^{k}(p-1)^{2}-4k(p+2)^{k-1} =\\
 &(k+3)(p+2)^{k-1}\left[ (1+\frac{1}{p+2})^{k-1}(p-1)^{2}(p+\frac{9-k}{3+k})-\frac{4k}{k+3}\right].
 \end{align*}
 If $k\geq 9$ there is nothing to prove. We assume $k=6, 7$ and $8$. In this case the only interesting situation is when $\frac{k-9}{3+k} \leq p \leq 0$. Then 
 \begin{align*}
&(1+\frac{1}{p+2})^{k-1}(p-1)^{2}(p+\frac{9-k}{3+k})-\frac{4k}{k+3} \leq \left(\frac{3}{2}\right)^{k-1}(p-1)^{2}(p+\frac{9-k}{3+k})-\frac{4k}{k+3} =\\
&\left(\frac{3}{2}\right)^{k-1}\left[p^{3}-\frac{3(k-1)}{k+3}p^{2}+\frac{3(k-5)}{k+3}p+ \frac{1}{k+3}\left( 9-k - 4k\left(\frac{2}{3}\right)^{k-1}\right)\right] \leq 0.
 \end{align*}
 The last inequality follows because the signs of the coefficients of the polynomial alternate, and $p\leq 0$. 
 
 \end{proof}

 \section{Conclusions}
 
 It follows from the previous section that if  $r_{1}> \ell r_{2}$ (where $\ell \approx 1.95$ is a positive solution of (\ref{shedegi})) then the heat loss of the big ball, i.e., $Q_{1}(d)$ is decreasing when $d \in [0,\varepsilon)$ where $\varepsilon>0$ is sufficiently small.
 Notice that when $d \to \infty$ we have $x \approx \frac{d}{r_{1}}$, and $y \approx \frac{d}{r_{2}}$ therefore  $\lim_{d \to \infty} f(d)=0$ and thus 
 \begin{align*}
 Q_{1}(\infty) = 4 \pi T_{0} r_{1} >  -4 \pi T_{0} r_{1}  \left(\frac{r_{2}(\gamma+\psi(\frac{r_{2}}{r_{1}+r_{2}}))}{r_{1}+r_{2}} \right)  = Q_{1}(0).
 \end{align*}
 The last inequality is justified because  $1>-x(\gamma+\psi(x))$ for $x \in (0,1)$ follows from the fact that $\psi(x+1)>-\gamma $ for $x \in (0,\infty)$. 
 So there exists a minimal value of $Q_{1}(d)$ on  $ [0, \infty)$, i.e., a nonzero distance when the heat loss of the big ball is minimal. 
 
\begin{figure}
\includegraphics[scale=0.4]{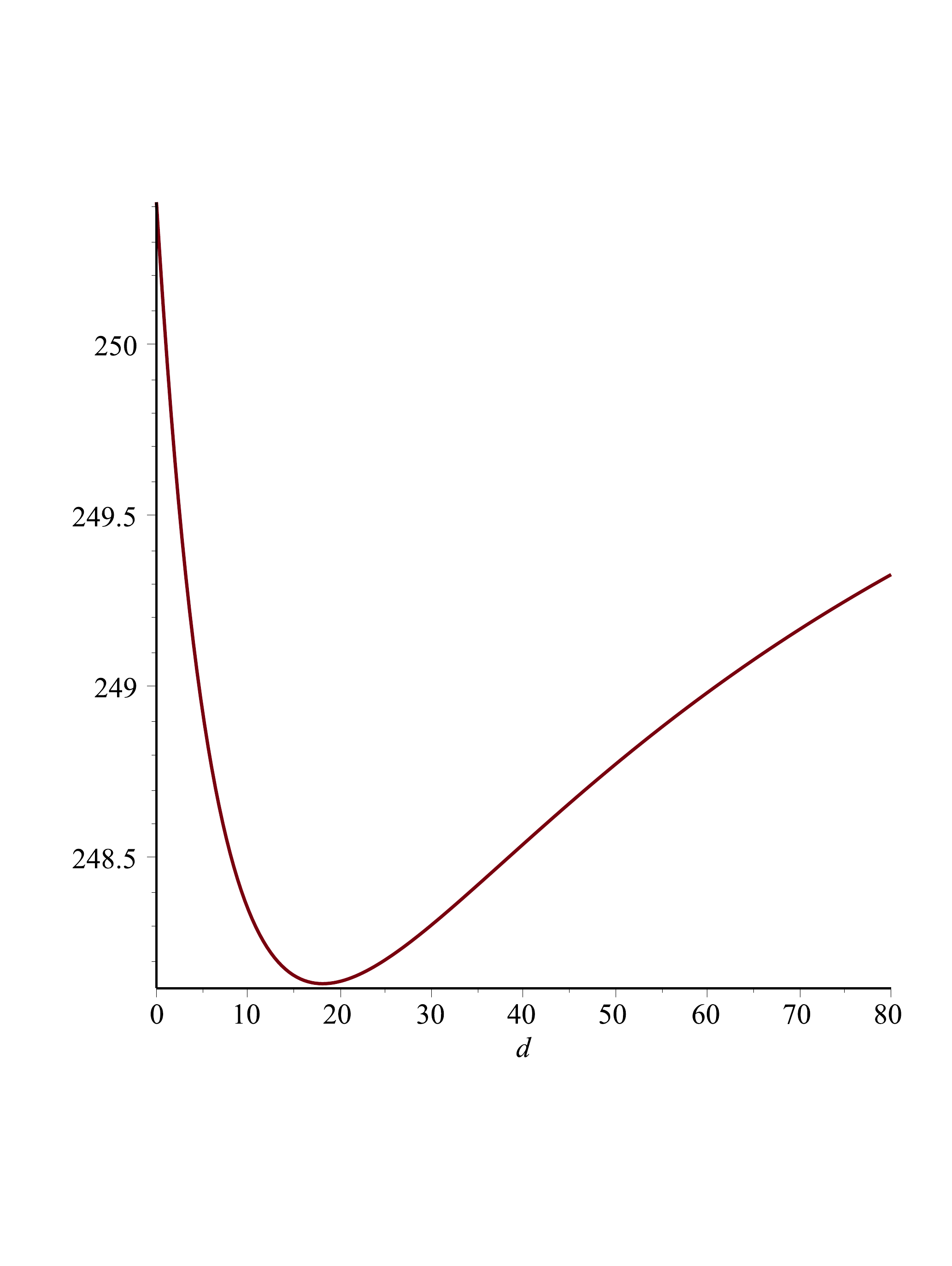}
\caption{The heat loss $Q_{1}(d)$ of the big ball. $T_{0}=1, r_{1}=20, r_{2}=1$, $0 \leq d \leq 80$}
\label{fig:basic}
\end{figure}

  The numerical computations show that, in fact  $Q_{1}(d)$ is decreasing on the interval $[0, c)$ and then it is increasing on $(c,\infty)$ where $c \approx r_{1}$. The heat loss of the small ball $Q_{2}(d)$ is always increasing for $d \geq 0$. Figure~\ref{fig:basic}  represents  the graph of  $Q_{1}(d)$  where $T_{0}=1, r_{1}=20, r_{2}=1$ and $0\leq t \leq 80$. 

 \subsection*{Acknowledgments}
 The author thanks Alexandre Eremenko, Benjamin Jaye and  Fedor Nazarov for helpful discussions and suggestions that led to simplifications of some proofs.

\end{document}